\documentclass[12pt]{article}
\usepackage{amsmath, amsfonts, amscd, latexsym, amsthm, amssymb}
\newtheorem{theor}{\hspace{1cm}{\sc Theorem}}

\newtheorem{sledst}[theor]{\hspace{1cm}Corollary}
\newtheorem{lemma}[theor]{\hspace{1cm}{\sc Lemma}}
\theoremstyle{definition}
\newtheorem{example}[theor]{\hspace{1cm}{\sc Example}}
\newtheorem{defin}[theor]{\hspace{1cm}{\sc Definition}}
\newtheorem*{rem}{\hspace{1cm}{\sc Remark}}
\newtheorem*{exa}{\hspace{1cm}{\sc Example}}

\newcommand{\Vol}{\mathop{\rm Vol}\nolimits}

\newcommand{\rk}{\mathop{\rm rk}\nolimits}

\newcommand{\conv}{\mathop{\rm conv}\nolimits}

\newcommand{\BP}{\mathop{\rm BP}\nolimits}
\newcommand{\MV}{\mathop{\rm MV}\nolimits}

\def\R{\mathbb R}
\def\N{\mathbb N}
\def\Z{\mathbb Z}

\def\C{\mathbb C}
\def\CC{({\mathbb C}\setminus\{0\})}

\begin{document}

\begin{center}
{\Large \textbf{Multiplicities of degenerations of matrices and
mixed volumes of Cayley polyhedra}}
\end{center}

\begin{center}
\textsc{A. Esterov}\footnote{This study was carried out within "The National Research University Higher School of Economics' Academic Fund Program in 2012-2013, research grant No.11-01-0125". Partially supported by RFBR,
grant 10-01-00678, MESRF, grant MK-6223.2012.1, and the Dynasty Foundation fellowship.}
\end{center}

\textbf{1. Introduction.}

The local version of D. Bernstein's formula \cite{bernst}
expresses the local degree of a germ of a proper analytic map in
terms of the Newton polyhedra of its components, provided that the
principal parts of its components are in general position (see
Theorem \ref{locbernst}). We generalize this formula as follows.

Let $A:\C^m\to\C^{n\times k}$ be a germ of an $(n\times
k)$-matrix with analytic entries, where $n\leqslant k$ (we denote
the space of all $(n\times k)$-matrices by $\C^{n\times k}$). If
$\rk A(0)<n$ and $\rk A(x)=n$ for all $x\ne 0$, then $m\leqslant
k-n+1$. Suppose that $m=k-n+1$ (in particular, if $n=1$, then
this means that $A:\C^m\to\C^m$ is a germ of a proper analytic
map). The intersection number $m(A)$ of the germ $A(\C^m)$ and
the set of all degenerate matrices in $\C^{n\times k}$ is well
defined, because the codimension of degenerate matrices in
$\C^{n\times k}$ equals $k-n+1$. In particular, if $n=1$, then
$m(A)$ equals the local degree of the map $A:\C^m\to\C^m$.

\begin{defin} \label{def1}
Let $A:\C^m\to\C^{n\times k}$ be a germ of an $(n\times
k)$-matrix with analytic entries, such that $m=k-n+1$, $\rk
A(0)<n$ and $\rk A(x)=n$ for all $x\ne 0$. Then the intersection
number $m(A)$ will be called \textit{the multiplicity (of degeneration)}
of the germ $A$.
\end{defin}
We recall the relation of this number to algebraic and topological invariants,
motivating our interest to it.

\textsc{Relation to Buchsbaum-Rim multiplicities.} In the notation of Definition
\ref{def1}, the multiplicity of the matrix $A$ is
equal to $\dim_{\C} \mathcal{O}_{\C^m,0}/\langle \mbox{maximal
minors of } A \rangle$, where $\mathcal{O}_{\C^m,0}$ is the ring
of germs of analytic functions on $\C^m$ near the origin. In particular,
it equals the Buchsbaum-Rim multiplicity of the submodule of $\mathcal{O}^n_{\C^m,0}$,
generated by the columns of $A$ (see, for example, Proposition 2.3 in \cite{gaff1}).

\textsc{Relation to characteristic classes.} Let $v_{i}$ be a
holomorphic section of a vector bundle $\mathcal{I}$ of rank $k$
on a smooth $(k-n+1)$-dimensional complex manifold $M$ for
$i=1,\ldots,n$. Suppose that there is a finite number of points
$x\in M$ such that the vectors $v_1(x),\ldots,v_{n}(x)$ are
linearly dependent. Denote the set of all such points by $X$.
Near each point $x\in X$, choosing a local basis $s_1,\ldots,s_k$
in the bundle $\mathcal{I}$, one can represent $v_i$ as a linear
combination $v_i=a_{i,1}s_1+\ldots+a_{i,k}s_k$, where $a_{i,j}$
are the entries of an $(n\times k)$-matrix $A:M\to\C^{n\times k}$
defined near $x$. Denote the multiplicity of $A$ by
$m_x$. Then the Chern number $c_{k-n+1}(\mathcal{I}_q)\cdot[M]$ is
equal to the sum of the multiplicities $m_x$ over all points $x\in
X$ (see, for example, \cite{gh}).

The aim of this paper is to present a formula for the
multiplicity of a matrix $A$ in terms of the Newton
polyhedra of the entries of $A$, provided that the principal
parts of the entries are in general position. In \cite{biv}, a
similar formula is given under the assumption that all the
entries from the same row of the matrix $A$ have the same Newton
polyhedron. \cite{E2} contains a general formula (see Theorem
\ref{th1st}), which is somewhat indirect in the sense that one
has to increase the dimension of polyhedra under consideration in
order to formulate the answer. The aim of this paper is to
simplify this answer combinatorially (see Theorem \ref{thmain}), so
that no higher-dimensional polyhedra are involved.

In Sections 2 and 3, we
present the formula for the multiplicity of a
matrix and the condition of general position for the principal
parts of the entries of a matrix, respectively. In Sections 5, 
this formula is deduced from Theorem
\ref{th1st}, which expresses the multiplicity of a matrix in terms of
the mixed volume of pairs of certain polyhedra (this notion is introduced in Section 4). This requires a formula
for the mixed volume of Cayley polyhedra (Theorem \ref{th2nd}, the proof given in Section 7),
which follows from the Oda equality $(A\cap\Z^n)+(B\cap\Z^n)=(A+B)\cap\Z^n$
for some class of bounded lattice polyhedra $A,B\subset\R^n$ (see
Section 6).

I am very grateful to the referee for many important remarks and ideas on how to improve the paper.


\textbf{2. Multiplicity in terms of Newton
polyhedra.}

A \textit{polyhedron} in $\R^n$ is the intersection of a finite
number of closed half-spaces. A \textit{face} of a polyhedron $A$
is the intersection of $A$ and the boundary of a closed
half-space, containing $A$. Note that the empty set is a face of
every polyhedron. \textit{The Minkowski sum} of sets $A$ and $B$
in $\R^n$ is the set $A+B=\{a+b\, |\, a\in A,\, b\in B\}$. Note
that $\varnothing+A=\varnothing$ for every $A$.

\begin{defin}
Let $B_i$ be a face of a polyhedron $\Delta_i\subset\R^m$ for
$i=1,\ldots,k$. The collection of faces $(B_1,\ldots,B_k)$ is
said to be \textit{compatible}, if the sum $B_1+\ldots+B_k$ is a
non-empty bounded face of the sum $\Delta_1+\ldots+\Delta_k$.
\end{defin}

Denote the positive orthants of $\R^m$ and $\Z^m$ by $\R^m_+$ and
$\Z^m_+$ respectively. For each point
$a=(a_1,\ldots,a_m)\in\Z^m$, denote the monomial $x_1^{a_1}\ldots
x_m^{a_m}$ by $x^a$.
\begin{defin}
\textit{The Newton polyhedron} $\Delta_f$ of a germ of an analytic
function $f=\sum_{a\in\Z^m_+}c_ax^a:\C^m\to\C$ is the convex hull
of the union $\bigcup\limits_{a\; |\; c_a\ne 0}{(a+\R^m_+)}$.
\end{defin}
\begin{defin}
\textit{The restriction} $f|_B$ of a germ
$f=\sum_{a\in\Z^m_+}c_ax^a$ to a bounded subset $B$ of the Newton
polyhedron $\Delta_f$ is the polynomial $\sum\limits_{a\in\Z^m\cap
B}c_ax^a$. The restriction of $f$ to the union of all bounded
faces of $\Delta_f$ is called \textit{the principal part} of $f$.
The restriction to the empty set equals zero by definition.
\end{defin}
The principal parts of the components of a map $f:\C^m\to\C^m$
form \textit{the principal part} of $f$, and the principal parts
of the entries of an $(n\times k)$-matrix $A:\C^m\to\C^{n\times
k}$ form \textit{the principal part} of $A$.

For a polyhedron $\Delta\subset\R^m_+$, denote the number of
integer lattice points in the difference $\R^m_+\setminus\Delta$
by $I(\Delta)$. Recall the local version of D. Bernstein's
formula \cite{bernst} (it can be deduced, for example from M.
Oka's formula \cite{oka}):
\begin{theor} \label{locbernst}
Let $f=(f_1,\ldots,f_m):\C^m\to\C^m$ be a germ of an analytic map
near the origin, and the differences
$\R^m_+\setminus\Delta_{f_i}$ are bounded. \newline 1) The local
degree of $f$ is greater than or equal to $$\sum_{0<p\leqslant m} (-1)^{m-p} \sum_{0<i_1<\ldots<i_p\leqslant
m}
I(\Delta_{f_{i_1}}+\ldots+\Delta_{f_{i_p}}),\eqno(*)$$
provided that $f$ is proper.
\newline 2) The germ $f$ is proper, and its local degree equals $(*)$, if and only if, for each compatible
collection of faces $B_1,\ldots,B_m$ of the polyhedra
$\Delta_1,\ldots,\Delta_m$, the system of polynomial equations
$f_1|_{B_1}=\ldots=f_m|_{B_m}=0$ has no roots in $\CC^m$.
\end{theor}
\begin{rem}
The principal parts, which satisfy the condition from part (2) of
this theorem, form a dense algebraic set in the space of
principal parts of maps with given Newton polyhedra of components.
\end{rem}
The main result of this paper is the following generalization of this fact to 
multiplicities of matrices.

\begin{defin} The \textit{tropical semiring} $P$ of polyhedra
is the set of all convex polyhedra in $\R^n$ (including the empty one)
with the additive operation $$A \vee B =\mbox{ convex hull of } A\cup B$$
and the Minkowski sum as  the multiplicative operation $$A + B = \{a+b\, |\, a\in A,\, b\in B\}.$$
\end{defin}
The name is justified by the fact that the support functions of $A \vee B$ and $A + B$ are equal to the
maximum and the sum of the support functions of $A$ and $B$ respectively. All the polyhedra $A$, satisfying the equation $A+\R^m_+=A$, form a subring $P_+\subset P$, and $\R^m_+$ is the unit in this  subring. In particular, whenever the sum of polyhedra $A_j\in P_+$ is taken over an empty set of indices $J=\varnothing$, we set $\sum_{j\in J} A_\alpha=\R^m_+$ by definition.
\begin{theor} \label{thmain}
Let $A=(a_{i,j}):\C^m\to\C^{n\times k}$ be a germ of an $(n\times
k)$-matrix with analytic entries, $m=k-n+1$, and the differences
$\R^m_+\setminus\Delta_{a_{i,j}}$ are bounded.  \newline 1) The
multiplicity of the matrix $A$ is greater than or equal to
$$\sum_{J\subset\{1,\ldots,k\}\atop b_1+\ldots+b_{n}=|J|}(-1)^{k-|J|}I\Bigl(\bigvee_{J_1\sqcup\ldots\sqcup J_n=J\atop |J_1|=b_1,\ldots,|J_n|=b_n}
\sum_{i=1,\ldots,n\atop j\in J_i}
\Delta_{a_{i,j}}\Bigr),\eqno(**)$$ provided that $\rk A(x)=n$ for
all $x\ne 0$. Here the first summation is taken over all non-empty $J\subset\{1,\ldots,k\}$ and all collections of 
non-negative integers $b_i$ that sum up to $|J|$, and 
$\bigvee$ is taken over all decompositions of $J$ into disjoint sets $J_i$ of size $b_i$.
\newline 2) We have $\rk
A(x)=n$ for all $x\ne 0$, and the multiplicity of
$A$ equals $(**)$, if and only if the principal part of $A$ is in general position 
in the sense of Definition \ref{nondegpp}.
\end{theor}

It is a purely combinatorial problem to deduce this fact from Theorem
\ref{th1st}, and it will be addressed in Section 5.

\begin{example} Theorem \ref{thmain} appears to be more convenient
than Theorem \ref{th1st} in many important special cases. For instance,
in the classical case
of homogeneous $a_{i,j}$, Theorem \ref{thmain} unlike Theorem \ref{th1st} gives a closed formula
for the multiplicity in terms of the degrees $d_{i,j}$ of the components $a_{i,j}$.  For $J\subset\{1,\ldots,k\}$ and a decomposition $|J|=b_1+\ldots+b_{n}$ into non-negative integers, introduce the number $$d^J_{b_1,\ldots,b_n}=\min_{J_1\sqcup\ldots\sqcup J_n=J\atop |J_1|=b_1,\ldots,|J_n|=b_n} \sum_{i=1,\ldots,n\atop j\in J_i} d_{i,j}.$$
\end{example}
\begin{sledst} In the setting of Theorem \ref{thmain}, assume that  the components $a_{i,j}$ are homogeneous polynomials of degree $d_{i,j}$.

1) The multiplicity of the matrix $A$ is greater or equal to
$$\sum_{J\subset\{1,\ldots,k\}\atop b_1+\ldots+b_{n}=|J|}(-1)^{k-|J|} {m+d^J_{b_1,\ldots,b_n}-1\choose m}.$$

2) The multiplicity is strictly greater than this number or is infinite, if and only if 
the entries are not in general position in the following sense: there exist integer numbers $\alpha_1,\ldots,\alpha_n$ and $\beta_1,\ldots,\beta_k$ and non-zero $x\in\C^m$ such that $d_{i,j}\geqslant\alpha_i+\beta_j$ for every $i$ and $j$, and the matrix of the entries $\delta_{d_{i,j}}^{\alpha_i+\beta_j}a_{i,j}(x)$ is effectively degenerate (as usual, $\delta_p^q$ is 1 if $p=q$ and 0 otherwise).
\end{sledst}

\begin{example} Note that, unlike in the complete intersection case $n=1$, the multiplicity of such a homogeneous matrix can be strictly greater than expected, but still finite. For example, if $(m,n,k)=(2,2,3)$, then the matrix $$ \left( \begin{array}{ccc}
 \;x+y & (x+y)^2+y^2 & x+y \\\
x+y & x+y & (x+y)^2+2y^2
\end{array} \right) $$
has multiplicity 6, which is strictly greater than the answer 3, given by Part 1 for a generic matrix of degree
$\left( \begin{array}{ccc}
 \;1 & 2 & 1 \\\
1 & 1 & 2
\end{array} \right). $
This is because the matrix above is not in general position (consider $\alpha_1=\alpha_2=1$, $\beta_1=\beta_2=\beta_3=0$ in the notation of Part 2). \end{example}

\begin{example} Let us expand the answer given by Theorem \ref{thmain} in the simplest case $(m,n,k)=(2,2,3)$.
Denote $\Delta_{a_{i,j}}$ by $\Delta_{i,j}$, then $(**)$ equals
$$I(\Delta_{1,1}+\Delta_{1,2}+\Delta_{1,3})+I\Bigl( (\Delta_{2,1}+\Delta_{1,2}+\Delta_{1,3})\vee
(\Delta_{1,1}+\Delta_{2,2}+\Delta_{1,3})\vee(\Delta_{1,1}+\Delta_{1,2}+\Delta_{2,3}) \Bigr)+$$
$$+I\Bigl( (\Delta_{1,1}+\Delta_{2,2}+\Delta_{2,3})\vee
(\Delta_{2,1}+\Delta_{1,2}+\Delta_{2,3})\vee(\Delta_{2,1}+\Delta_{2,2}+\Delta_{1,3}) \Bigr)+
I(\Delta_{2,1}+\Delta_{2,2}+\Delta_{2,3})-$$
$$-I(\Delta_{1,1}+\Delta_{1,2})-I(\Delta_{1,1}+\Delta_{1,3})-I(\Delta_{1,2}+\Delta_{1,3})-I(\Delta_{2,1}+\Delta_{2,2})
-I(\Delta_{2,1}+\Delta_{2,3})-I(\Delta_{2,2}+\Delta_{2,3})-$$
$$-I( (\Delta_{1,1}+\Delta_{2,2})\vee(\Delta_{1,2}+\Delta_{2,1}))
-I( (\Delta_{1,1}+\Delta_{2,3})\vee(\Delta_{1,3}+\Delta_{2,1}))
-I( (\Delta_{1,3}+\Delta_{2,2})\vee(\Delta_{1,2}+\Delta_{2,3}))-$$
$$+I(\Delta_{1,1})+I(\Delta_{1,2})+I(\Delta_{1,3})+I(\Delta_{2,1})+I(\Delta_{2,2})+I(\Delta_{2,3}).$$ \end{example}
\begin{example} 
If $\Delta_{i,j}=\Delta_i$ does not depend on the column $j$, then the answer, given by Theorems \ref{thmain} and \ref{th1st},
admits a much simpler form $\sum_{1\leqslant i_1\leqslant\ldots\leqslant i_m\leqslant k} \MV(\Delta_{i_1},\ldots,\Delta_{i_m})$. If $\Delta_{i,j}=\Delta_j$ does not depend on the row $i$, then the answer, given by Theorems \ref{thmain} and \ref{th1st},
admits a much simpler form $\sum_{1\leqslant j_1<\ldots<j_m\leqslant k} \MV(\Delta_{j_1},\ldots,\Delta_{j_m})$. 
Both of these facts can be easily deduced from Theorem \ref{th1st} (see \cite{E3} and \cite{E4} for details). 
The latter one was discovered earlier in a much more general setting by Bivi\`{a}-Ausina (\cite{biv}).

For example, if the germs $a_{i1}\in\langle x^2,y\rangle,\, a_{i2}\in\langle x,y^3\rangle,\, a_{i3}\in\langle x^2,y^3\rangle$, $i=1,2$, are in general position, then the multiplicity of $A$ equals $4I(\Delta_1+\Delta_2+\Delta_3)-3I(\Delta_1+\Delta_2)-3I(\Delta_1+\Delta_3)-3I(\Delta_2+\Delta_2)
+2I(\Delta_1)+2I(\Delta_2)+2I(\Delta_3)=4\cdot 16-3\cdot (6+9+11)+2\cdot (2+3+5)=6$ according to Theorem \ref{thmain}
and $\MV(\Delta_1,\Delta_2)+\MV(\Delta_1,\Delta_3)+\MV(\Delta_2,\Delta_3)=1+2+3=6$ according to \cite{biv}.
\end{example}

\textbf{3. General position of principal parts of matrices.}

By convention, each polyhedron has the empty face. In particular,
some faces $B_{i,j}$ in the following definition may be empty.

\begin{defin}\label{matrcomdef} Let $B_{i,j}$ be a bounded face of a polyhedron
$\Delta_{i,j}\subset\R^m$ for $i=1,\ldots,n$, $j=1,\ldots,k$. The
collection of faces $B_{i,j}$ is said to be
\textit{matrix-compatible}, if there exist vectors
$c_1,\ldots,c_n\in\Z^m$ and compatible faces $B_1,\ldots,B_k$ of
the convex hulls
$\bigvee_i(\Delta_{i,1}+c_i),\ldots,\bigvee_i(\Delta_{i,k}+c_i)$,
such that $B_{i,j}=(B_j-c_i)\cap\Delta_{i,j}$ for each
$i=1,\ldots,n$, $j=1,\ldots,k$.
\end{defin}

\begin{example} Let $\Delta_{i,j}\subset\R^1$ be the rays
$$ \left( \begin{array}{ccc}
 \;[1,\infty) & [1,\infty) & [1,\infty) \\\
[1,\infty) & [2,\infty) & [2,\infty) \\\
[1,\infty) & [2,\infty) & [2,\infty)
\end{array} \right), $$
then every face $B_{i,j}$ is either the origin of $\Delta_{i,j}$
(denoted by $*$), or empty (denoted by $\varnothing$). In this
case, the matrix-compatible collection of faces are
$$\mathcal{B}_1=\left( \begin{array}{ccc}
* & * & * \\
* & \varnothing & \varnothing \\
* & \varnothing & \varnothing \end{array} \right),\quad \mathcal{B}_2=\left( \begin{array}{ccc}
\varnothing & * & * \\
* & * & * \\
* & * & * \end{array} \right), $$
and 11 more collections with fewer non-empty faces.\end{example}
\begin{defin} \label{nondegmat}
A matrix $M\in\C^{n\times k},\; n\leqslant k$, is said to be
\textit{effectively non-degenerate}, if $(t_1,\ldots,t_n)\cdot
M\ne(0,\ldots,0)$ for all $(t_1,\ldots,t_n)\in\CC^n$.
\end{defin}
\begin{example} The complex matrix
$$\left( \begin{array}{ccc}
a & b & c \\
d & 0 & 0 \\
e & 0 & 0 \end{array} \right)$$
is effectively non-degenerate if and only if $b=c=0$ (although it is degenerate for all complex numbers $a,b,c,d,e$).

For an $(n\times k)$-matrix $A$ with analytic entries
$a_{i,j}:\C^m\to\C$ and a collection $\mathcal{B}$ of faces
$B_{i,j}$ of the Newton polyhedra $\Delta_{a_{i,j}}$, we denote the
matrix with entries $a_{i,j}|_{B_{i,j}}$ by $A|_{\mathcal{B}}$.\end{example}

\begin{defin} \label{nondegpp}
The principal part of an $(n\times k)$-matrix $A$ with analytic
entries $a_{i,j}:\C^m\to\C$ is said to be
\textit{in general position}, if, for each matrix-compatible
collection $\mathcal{B}$ of faces of the Newton polyhedra
$\Delta_{a_{i,j}}$ and for each $x\in\CC^m$, the matrix
$A|_{\mathcal{B}}(x)$ is effectively non-degenerate.
\end{defin}
\begin{rem}
Principal parts in general position form a dense algebraic set in the
space of principal parts of matrices with given Newton polyhedra
of entries. However, this is not true, if we replace the effective 
non-degeneracy of matrices with the conventional one in 
Definition \ref{nondegpp}. For instance, if $(m,n,k)=(1,3,3)$,
and the Newton polyhedra $\Delta_{a_{i,j}}$ are as in the example 
to Definition \ref{matrcomdef}, then the only non-trivial condition,
imposed by Definition \ref{nondegpp}, corresponds to the second
matrix-compatible collection of faces shown in the example:
$$\det(A|_{\mathcal{B}_2})=\det\left( \begin{array}{ccc}
0 & a^0_{0,1} & a^0_{0,2} \\
a^0_{1,0} & a^0_{1,1} & a^0_{1,2} \\
a^0_{2,0} & a^0_{2,1} & a^0_{2,2} \end{array} \right)\ne 0,$$
where $a^0_{i,j}$ is the leading coefficient of the series $a_{i,j}$.
However, if we replace effective nondegeneracy with nondegeneracy
in Definition \ref{nondegpp}, then no matrix $A$ will satisfy it,
because the matrix $A|_{\mathcal{B}_1}$ is always degenerate
(see the example to Definition \ref{nondegmat}).

It would be thus interesting to describe
a collection of minors of the matrices $A|_{\mathcal{B}}$,
such that  \newline 1) If the
principal part of $A$ is in general position, then these minors vanish. \newline 2) The principal
parts for which these minors vanish form a (closed algebraic)
set of positive codimension in the space of all principal
parts of matrices with given Newton polyhedra of entries.

This reduces to the following problem: given $K\subset\N^2$, assume that $a_{i,j}$
are independent variables for $(i,j)\in K$, and the entries of the matrix $A$ 
are $a_{i,j}$ for $(i,j)\in K$ and equal 0 for $(i,j)\notin K$. Find a collection 
of minors $\mathcal{A}$
of the matrix $A$, such that \newline 1) If $A$ is effectively nondegenerate, then $\mathcal{A}=0$. \newline 2) We have
$\mathcal{A}\ne 0$ for generic $a_{i,j}, (i,j)\in K$.
\end{rem}

\textbf{4. Mixed volumes of pairs of polyhedra.}

\begin{defin}
Polyhedra $\Delta_1$ and $\Delta_2$ in $\R^n$ are said to be
\textit{parallel} if $a+\Delta_1\subseteq\Delta_1 \Leftrightarrow
a+\Delta_2\subseteq\Delta_2$ for every point $a\in\R^n$.
\end{defin}

\begin{defin} \label{defmixvol}  (\cite{E2}, \cite{E3}) 1) A pair of polyhedra $\Delta_1,\Delta_2$ in $\R^n$ is
called \textit{bounded} if both $\Delta_1\setminus\Delta_2$ and
$\Delta_2\setminus\Delta_1$ are bounded. The set of all bounded
pairs of polyhedra parallel to a given convex cone $C\subset\R^n$
is denoted by $\BP_C$.

2) \textit{The Minkowski sum} $(\Delta_1, \Delta_2) + (\Gamma_1,
\Gamma_2)$ of two pairs from $\BP_C$ is the pair
$(\Delta_1+\Gamma_1, \Delta_2+\Gamma_2)\in \BP_C$.

3) \textit{The volume} $\Vol(\Delta_1,\Delta_2)$ of a bounded pair
$(\Delta_1,\Delta_2)$ is the difference
$\Vol(\Delta_1\setminus\Delta_2)-\Vol(\Delta_2\setminus\Delta_1)$.

4) \textit{The mixed volume} is the symmetric multilinear (with
respect to Minkowski summation) function
$\MV:\underbrace{\BP_C\times\ldots\times\BP_C}_n\to\R$ such that
$\MV(A,\ldots,A)=\Vol(A)$ for every pair $A\in\BP_C$.
\end{defin}

There exists a unique such function $\MV$ (see \cite{E3}, Section 4, Lemma 3
for existance, uniqueness and all other basic
facts about the mixed volume of pairs, mentioned below). Recall that a polyhedron
is said to be \textit{lattice} if its vertices are integer lattice
points. The mixed volume of pairs of
$n$-dimensional lattice polyhedra is a rational number with
denominator $n!$. 

\begin{exa} If $C$ consists of one point, then $\BP_C$ consists of
pairs of bounded polyhedra, and $$\MV\bigl( (\Delta_1, \Gamma_1),
\ldots, (\Delta_n, \Gamma_n)
\bigr)=\MV(\Delta_1,\ldots,\Delta_n)-\MV(\Gamma_1,\ldots,\Gamma_n),$$
where $\MV$ in the right hand side is the classical mixed volume
of bounded polyhedra. If $C$ is not bounded, then both terms in
the right hand side are infinite, but "their difference makes
sense".
\end{exa}

One can use the following formula to express the mixed volume of
pairs in terms of mixed volumes of polyhedra (\cite{E3}, Section 4, Lemma 3).

\begin{lemma}
For bounded pairs $(\Delta_i,\Gamma_i)\in\BP_C,\, i=1,\ldots,n$,
let $H\subset\R^n$ be a half-space such that $C\cap H$ is bounded
and $\Delta_i\setminus H=\Gamma_i\setminus H$. Then
$$\MV\bigl( (\Delta_1, \Gamma_1), \ldots, (\Delta_n, \Gamma_n)
\bigr)=\MV(\Delta_1\cap H,\ldots,\Delta_n\cap H)-\MV(\Gamma_1\cap
H,\ldots,\Gamma_n\cap H),$$ where $\MV$ in the right hand side is
the classical mixed volume of bounded polyhedra.
\end{lemma}

For a
bounded pair of (closed) polyhedra $(\Delta,\Gamma)\in\BP_C$,
define $I(\Delta,\Gamma)$ as the number of integer lattice points
in the difference $\Delta\setminus\Gamma$ minus the number of
integer lattice points in the difference $\Gamma\setminus\Delta$.

\begin{lemma} \label{volint}
For bounded pairs of lattice polyhedra $A_i\in\BP_C$, we have
$$n!\MV(A_1,\ldots,A_n)=\sum_{0<p\leqslant m} (-1)^{n-p} \sum_{0<i_1<\ldots<i_p\leqslant n}
I(A_{i_1}+\ldots+A_{i_p}).$$
\end{lemma}
\textsc{Proof.} For the classical mixed volume of bounded polyhedra,
this equality is well known (see, for example, \cite{khovpp}). The general
case can be deduced to the case of bounded polyhedra by the previous lemma. $\quad\Box$

\textbf{5. Proof of Theorem \ref{thmain}.}

The following theorem is a special case of Theorem 5 from
\cite{E3}.

\begin{defin}
For polyhedra $\Delta_1,\ldots,\Delta_n\subset\R^m$, define
the \textit{Cayley polyhedron} $\Delta_1*\ldots*\Delta_n$ 
as the convex hull of the union
$$\bigcup_i \{b_i\}\times\Delta_i\subset\R^{n-1}\oplus\R^m,$$
where $b_1,\ldots,b_n$ are the points $(1,0,\ldots,0),
(0,1,\dots,0), \ldots, (0,0,\dots,1)$ and $(0,0,\dots,0)$ in
$\R^{n-1}$. Denote $\R^m_+*\ldots*\R^m_+$ by $D$.
\end{defin}

For germs of analytic functions $a_1,\ldots,a_n$ on $\C^m$ near
the origin, denote the sum $t_1a_1+\ldots+t_{n-1}a_{n-1}+a_n$ by
$a_1*\ldots*a_n$, where $t_1,\ldots,t_{n-1}$ are the standard
coordinates on $\C^{n-1}$.

\begin{theor} (\cite{E2}, \cite{E3}, \cite{E4}) \label{th1st}
Let $A$ be an $(n\times k)$-matrix with entries
$a_{i,j}:\C^m\to\C$ which are germs of analytic functions near
the origin. Suppose that the Newton polyhedra $\Delta_{i,j}$ of
the germs $a_{i,j}$ intersect all coordinate axes in $\R^m$.
\newline 1) The multiplicity of $A$ is greater than or equal to
$$(m+n-1)!\MV\bigl((D,\Delta_{1,1}*\ldots*\Delta_{n,1}),\ldots,(D,\Delta_{1,k}*\ldots*\Delta_{n,k})\bigr).\eqno(***)$$
2) We have $\rk
A(x)=n$ for all $x\ne 0$, and the multiplicity of $A$
equals $(***)$, if and only if, for each compatible collection of faces $B_1,\ldots,B_k$
of the polyhedra ${\Delta_{1,1}*\ldots*\Delta_{n,1}},$ $\ldots,\,
{\Delta_{1,k}*\ldots*\Delta_{n,k}}$, the polynomials
${(a_{1,1}*\ldots*a_{n,1})|_{B_1}},\,$ $\ldots,\,$
${(a_{1,k}*\ldots*a_{n,k})|_{B_k}}$ have no common zeroes in
$\CC^{n-1}\times\CC^m$.
\end{theor}
The ``only if'' part of (2) is actually proved in \cite{E3}, but is
explicitly formulated and discussed only in \cite{E4}, Theorem 1.21.

Recall that $|S\cap\Z^m|$ is denoted by $I(S)$ for a bounded set $S\in\R^m$. 
If the symmetric difference of
(closed) lattice polyhedra $\Gamma$ and $\Delta$ in $\R^m$ is
bounded, denote the difference
$I(\Gamma\setminus\Delta)-I(\Delta\setminus\Gamma)$ by
$I(\Gamma,\Delta)$. 
For pairs
of polyhedra $(\Gamma_i,\Delta_i)$ in $\R^m$, denote the pair
$(\bigvee_i\Gamma_i,\bigvee_i\Delta_i)$ by
$\bigvee_i(\Gamma_i,\Delta_i)$ and  the
pair $(\Gamma_1*\ldots*\Gamma_n,\Delta_1*\ldots*\Delta_n)$ 
by $(\Gamma_1,\Delta_1)*\ldots*(\Gamma_n,\Delta_n)$. \begin{theor} \label{th2nd} If
$B_{i,j}$, $i=1,\ldots,n$, $j=1,\ldots,k$, are bounded lattice polyhedra in $\R^m$ or 
pairs of lattice polyhedra in $\BP_C$, and $m=k-n+1$, then the mixed volume of
$B_{1,j}*\ldots*B_{n,j}$, $j=1,\ldots,k$, equals
$$\frac{1}{k!}\sum_{J\subset\{1,\ldots,k\}\atop b_1+\ldots+b_{n}=|J|}(-1)^{k-|J|}I\Bigl(\bigvee_{J_1\sqcup\ldots\sqcup J_n=J\atop |J_1|=b_1,\ldots,|J_n|=b_n}
\sum_{i=1,\ldots,n\atop j\in J_i} B_{i,j}\Bigr).$$\end{theor} 
Note that some of $B_{i,j}$ may be empty. The proof is given in Section 7.
Theorem \ref{thmain} follows from Theorems \ref{th1st} and \ref{th2nd}
(one can easily check that the condition of general position in Theorem \ref{th1st}(2)
coincides with the one given by Definition \ref{nondegpp}).

\textbf{6. Fans and lattice points of polyhedra.}

Here we prove the equality
$$(A\cap\Z^q)+(B\cap\Z^q)=(A+B)\cap\Z^q$$ for some class of
bounded lattice polyhedra $A,B\subset\R^q$ (see \cite{oda1}
for a conjecture in the general case). 

\begin{defin} \label{defcone} \textit{A (rational) cone} in $\R^q$ generated by (rational) vectors
$v_1,\ldots,v_m$ is the set of all linear combinations of
$v_1,\ldots,v_m$ with positive coefficients.
\end{defin}

Note that, according to this definition, a cone is not a closed set unless it is a vector subspace of $\R^q$,
and is not an open set unless it is $q$-dimensional.

\begin{defin}
A collection of rational cones $C_1,\ldots,C_p$ in $\R^q$ is said
to be $\Z$-\textit{transversal}, if $\sum\dim C_i=q$ and the set
$\Z^q\cap\bigcup_i {C_i}$ generates the lattice $\Z^q$.
\end{defin}


\begin{defin}
\textit{A (rational) fan} $\Phi$ in $\R^q$ is a non-empty finite
set of nonoverlapping (rational) cones in $\R^q$ such that
\newline 1) Each face of each cone from $\Phi$ is in $\Phi$,
\newline 2) Each cone from $\Phi$ is a face of a $q$-dimensional
cone from $\Phi$.
\end{defin}

\begin{defin}
A collection of fans $\Phi_1,\ldots,\Phi_p$ in $\R^q$ is said to
be $\Z$-\textit{transversal with respect to shifts} $c_1\in\R^q,
\ldots, c_p\in\R^q$, if each collection of cones $C_1\in\Phi_1,
\ldots, C_p\in\Phi_p$, such that the intersection
$(C_1+c_1)\cap\ldots\cap (C_p+c_p)$ consists of one point, is
$\Z$-transversal.
\end{defin}

\begin{defin} \label{dualfan}
\textit{The dual cone} of a face $B$ of a polyhedron
$A\subset\R^q$ is the set of all covectors $\gamma\in(\R^q)^*$
such that $\{a\in A\, |\, \gamma(a)=\min\gamma(A)\}=B$.
\textit{The dual fan} of a polyhedron is the set of dual cones of
all its faces.
\end{defin}

\begin{theor} \label{transvpoly}
If the dual fans of bounded lattice polyhedra
$A_1,\ldots,A_p\subset\R^q$ are $\Z$-transversal with respect to
some shifts $c_1\in(\R^q)^*, \ldots, c_p\in(\R^q)^*$ and
$\dim(A_1+\ldots+A_p)=q$, then
$$(A_1\cap\Z^q)+\ldots+(A_p\cap\Z^q)=(A_1+\ldots+A_p)\cap\Z^q.$$
\end{theor}

\textsc{Proof.} Consider covectors $c_1\in(\R^q)^*, \ldots,
c_p\in(\R^q)^*$ as linear functions on the polyhedra
$A_1\subset\R^q,\ldots,A_p\subset\R^q$ respectively, and denote
their graphs in $\R^q\oplus\R^1$ by $\Gamma_1,\ldots,\Gamma_p$.
Denote the projection $\R^q\oplus\R^1\to\R^q$ by $\pi$, and
denote the ray $\{ (0,\ldots,0,t)\; |\;
t<0\}\subset\R^q\oplus\R^1$ by $L_-$.

Each bounded $q$-dimensional face $B$ of the sum
$\Gamma_1+\ldots+\Gamma_p+L_-$ is the sum of some faces
$B_1,\ldots,B_p$ of polyhedra $\Gamma_1+L_-,\ldots,\Gamma_p+L_-$.
$\Z$-transversality with respect to shifts $c_1\in(\R^q)^*,
\ldots, c_p\in(\R^q)^*$ implies that
$$\bigl(\pi(B_1)\cap\Z^q\bigr)+\ldots+\bigl(\pi(B_p)\cap\Z^q\bigr)=\pi(B_1+\ldots+B_p)\cap\Z^q.$$
Since the projections of bounded $q$-dimensional faces of the sum
$\Gamma_1+\ldots+\Gamma_p+L_-$ cover the sum $A_1+\ldots+A_p$, it
satisfies the same equality:
$$(A_1\cap\Z^q)+\ldots+(A_p\cap\Z^q)=(A_1+\ldots+A_p)\cap\Z^q.\eqno\Box$$

\begin{sledst} \label{minimax1}
Let $S\subset\R^q$ be the standard $q$-dimensional simplex, let
$l_1,\ldots,l_p$ be linear functions on $S$ with graphs
$\Gamma_1,\ldots,\Gamma_p$, and let $l$ be the maximal
piecewise-linear function on $pS$, such that its graph $\Gamma$
is contained in the sum $\Gamma_1+\ldots+\Gamma_p$. Then, for
each integer lattice point $a\in pS$, the value $l(a)$ equals the
maximum of sums $l_1(c_1)+\ldots+l_p(c_p)$, where
$(c_1,\ldots,c_p)$ runs over all $p$-tuples of vertices of $S$
such that $c_1+\ldots+c_p=a$.
\end{sledst}

\textsc{Proof.} Denote the projection $\R^q\oplus\R^1\to\R^q$ by
$\pi$. A $q$-dimensional face $B$ of $\Gamma$, which contains the
point $\bigl(a,l(a)\bigr)\in\R^q\oplus\R^1$, can be represented
as a sum of faces $B_i$ of simplices $\Gamma_i$. Since
$\pi(B_1),\ldots,\pi(B_p)$ are faces of the standard simplex,
their dual fans are $\Z$-transversal with respect to a generic
collection of shifts, and, by Theorem \ref{transvpoly},
$$\bigl(\pi(B_1)\cap\Z^q\bigr)+\ldots+\bigl(\pi(B_p)\cap\Z^q\bigr)=\pi(B)\cap\Z^q.$$
In particular, $a=c_1+\ldots+c_p$ for some integer lattice points
$c_i\in\pi(B_i)$, which implies $l(a)=l_1(c_1)+\ldots+l_p(c_p)$.
$\Box$

\begin{rem}
In particular, if the functions $l_1,\ldots,l_p$ are in general
position, then all $C_{p+q}^{q}$ integer lattice points in the
simplex $pS$ are projections of vertices of $\Gamma$. Translating
this into the tropical language, one can prove again the following
well-known fact: $p$ generic tropical hyperplanes in the space
$\R^q$ subdivide it into $C_{p+q}^{q}$ pieces.
\end{rem}

\begin{example}
If $S$ in the formulation of Corollary \ref{minimax1} is not the
standard simplex, then the statement is not always true. For
example, consider $$S=\conv\bigl\{(1,1), (1,-1), (-1,1),
(-1,-1)\bigr\},\; l_1(x,y)=x+y,\; l_2(x,y)=x-y,\; a=(1,0).$$

If, in addition, we allow functions $l_j$ to be concave piecewise
linear with integer domains of linearity, then the statement is
not true unless $S$ is the standard simplex. That is why we cannot
use computations below to simplify the formula in the statement of
Theorem 5 from \cite{E3} in general.
\end{example}

\textbf{7. Proof of Theorem \ref{th2nd}.}

Rewriting the mixed volume of the pairs $B_{1,i}*\ldots*B_{n,i},\, i=1,\ldots,k,$ as
$$\sum_{0<p\leqslant m} (-1)^{n-p} \sum_{0<i_1<\ldots<i_p\leqslant n}
I\bigl((B_{1,i_1}*\ldots*B_{n,i_1})+\ldots+(B_{1,i_p}*\ldots*B_{n,i_p})\bigr)$$
by Lemma \ref{volint}, and applying the following Lemma \ref{intsum}
to every term in this sum, we obtain the statement of Theorem \ref{th2nd}.



\begin{lemma} \label{intsum}
For bounded pairs of polyhedra
$A_{i,j}=(\Delta_{i,j},\Phi_{i,j})\in\BP_C$, $i=1,\ldots,n$,
$j=1,\ldots,p$,
$$I\bigl((A_{1,1}*\ldots*A_{n,1})+\ldots+(A_{1,p}*\ldots*A_{n,p})\bigr)=$$
$$=\sum_{a_1+\ldots+a_n=p\atop a_1\geqslant 0,\ldots, a_n\geqslant 0}I\Bigl(\bigvee_{J_1\sqcup\ldots\sqcup J_n=\{1,\ldots,p\}\atop |J_1|=a_1,\ldots,|J_n|=a_n}
\sum_{i=1,\ldots,n\atop j\in J_i} \Delta_{i,j},\;
\bigvee_{J_1\sqcup\ldots\sqcup J_n=\{1,\ldots,p\}\atop
|J_1|=a_1,\ldots,|J_n|=a_n} \sum_{i=1,\ldots,n\atop j\in J_i}
\Phi_{i,j}\Bigr).$$
\end{lemma}

\textsc{Proof.} Every integer lattice point, participating in
the left hand side, is contained in the plane
$\{(a_1,\ldots,a_{n-1})\}\times\R^m\subset\R^{n-1}\oplus\R^m$ for
some non-negative integer numbers $a_1,\ldots,a_n$, which sum up
to $p$. Thus, it is enough to describe the intersection of the
pair
$\bigl((A_{1,1}*\ldots*A_{n,1})+\ldots+(A_{1,p}*\ldots*A_{n,p})\bigr)$
with each of these planes, using the following fact. $\Box$

\begin{lemma} Suppose that polyhedra $\Delta_{i,j}\subset\R^m$ are parallel to each other for $i=1,\ldots,n$,
$j=1,\ldots,p$. Then, for each $n$-tuple of non-negative integer
numbers $a_1,\ldots,a_n$ which sum up to $p$,
$$\Bigl(\{(a_1,\ldots,a_{n-1})\}\times\R^m\Bigr)\, \cap\,
\bigl((\Delta_{1,1}*\ldots*\Delta_{n,1})+\ldots+(\Delta_{1,p}*\ldots*\Delta_{n,p})\bigr)=$$
$$=\{(a_1,\ldots,a_{n-1})\}\times\Bigl(\bigvee_{J_1\sqcup\ldots\sqcup J_n=\{1,\ldots,p\}\atop |J_1|=a_1,\ldots,|J_n|=a_n}
\sum_{i=1,\ldots,n\atop j\in J_i}
\Delta_{i,j}\Bigr)\subset\R^{n-1}\oplus\R^m.$$
\end{lemma}

\textsc{Proof.} For each hyperplane $L\subset\R^m$, denote the
projection $\R^{n-1}\oplus\R^m\to\R^{n-1}\oplus\R\, $ along $\,
\{0\}\oplus L$ by $\pi_L$. It is enough to prove that the images
of the left hand side and the right hand side under $\pi_L$
coincide for each $L$. To prove it, apply Corollary
\ref{minimax1}, setting $q$ to $n-1$, $a$ to
$(a_1,\ldots,a_{n-1})$, and $\Gamma_j$ to the maximal bounded
face of the projection
$\pi_L\bigl(\Delta_{1,j}*\ldots*\Delta_{n,j}\bigr)$ for every
$j=1,\ldots,p$. $\Box$

\noindent(A. Esterov) \textsc{National Research University Higher School of Economics. \newline Faculty of Mathematics NRU HSE, 7 Vavilova 117312 Moscow, Russia.}


\begin{thebibliography}{9999}

\bibitem[Ber]{bernst}
D. N. Bernstein; The number of roots of a system of equations;
Functional Anal. Appl. 9 (1975), no. 3, 183--185.

\bibitem[Biv]{biv}
C. Bivi\`{a}-Ausina; The integral closure of modules,
Buchsbaum-Rim multiplicities and Newton polyhedra;
J.~London~Math.~Soci. (2) 69 (2004), 407--427.

\bibitem[EG]{smg}
W. Ebeling, S. Gusein-Zade; Indices of vector fields or 1-forms
and characteristic numbers; Bull. London Math. Soc. 37 (2005),
no. 5, 747--754.

\bibitem[G]{gaff1}
T. Gaffney; Multiplicities and equisingularity of ICIS germs; Invent. Math. 123 (1996), 209--220.

\bibitem[O97]{oda1}
T. Oda; Problems on Minkowski sums of convex lattice polytopes;
arXiv:0812.1418 (1997).

\bibitem[O90]{oka}
M. Oka; Principal zeta-function of non-degenerate complete
intersection singularity; J.~Fac.~Sci.~Univ.~Tokyo 37 (1990),
11--32.

\bibitem[S05]{suwa}
T. Suwa; Residues of Chern classes on singular varieties;
Singularit\'es Franco-Japonaises, 265--285, S\'emin. Congr., 10,
Soc. Math. France, Paris, 2005.

\bibitem[Kh]{khovpp}
A. G. Khovanskii; Newton polyhedra, and the genus of complete
intersections; Functional Anal. Appl. 12 (1978), no. 1, 38--46.

\bibitem[GH]{gh}
P. Griffiths, J. Harris; Principles of algebraic geometry; John
Wiley \& Sons, New York, 1978.

\bibitem[E05]{E2}
A. Esterov; Indices of 1-forms, resultants, and Newton polyhedra;
Russian Math. Surveys 60 (2005), no. 2, 352--353.

\bibitem[E06]{E3}
A. Esterov; Indices of 1-forms, intersection indices, and Newton
polyhedra; Sb. Math., 197 (2006), no. 7, 1085--1108

\bibitem[E09]{E4}
A. Esterov; Determinantal singularities and Newton
polyhedra; arXiv:0906.5097

\end{thebibliography}
\end{document}